\documentclass[11pt]{amsart}
\usepackage[latin1]{inputenc}
\usepackage{amsfonts}
\usepackage[toc,page,title,titletoc,header]{appendix}
\usepackage{enumerate,enumitem} 
\usepackage{graphicx,psfrag,epsfig,multirow,caption,subcaption}
\usepackage{lmodern,microtype,amssymb,amsmath,amscd,amsthm,amssymb,verbatim,setspace,mathtools}
\numberwithin{equation}{section}
\usepackage{mathrsfs,wrapfig}
\usepackage{indentfirst}
\usepackage{extarrows}
\usepackage{float}
\usepackage{xcolor}
\usepackage[colorlinks=true, backref=page]{hyperref}
\usepackage{marginnote}
\usepackage[margin=1.08in]{geometry}

\numberwithin{equation}{section}
\newtheorem{theorem}{Theorem}[section]
\newtheorem{proposition}[theorem]{Proposition}
\newtheorem{lemma}[theorem]{Lemma}
\newtheorem{corollary}[theorem]{Corollary}

\newtheorem{Claim}[theorem]{Claim}
\theoremstyle{definition}
\theoremstyle{remark}\newtheorem{remark}[theorem]{Remark}

\newcommand{\C}{\mathbb C}

\newcommand{\R}{\mathbb R}
\newcommand{\Sph}{\mathbb S}

\newcommand{\Z}{\mathbb Z}

\newcommand{\bx}{\mathbf{x}}
\newcommand{\by}{\mathbf{y}}

\newcommand{\cH}{\mathcal{H}}

\newcommand{\cL}{\mathcal{L}}

\newcommand{\cP}{\mathcal{P}}

\newcommand{\mfk}{\mathfrak}

\newcommand{\eps}{\varepsilon}

\newcommand{\al}{\alpha}
\newcommand{\Euc}{\mathrm{euc}}

\DeclareMathOperator{\vol}{Vol}

\DeclareMathOperator{\spt}{spt}

\DeclareMathOperator{\diag}{diag}

\DeclareMathOperator{\dist}{dist}

\DeclareMathOperator{\Tan}{Tan}

\DeclareMathOperator{\loc}{loc}
\DeclareMathOperator{\tr}{tr}

\DeclareMathOperator{\Reg}{reg}
\DeclareMathOperator{\Sing}{sing}

\DeclareMathOperator{\re}{Re}
\DeclareMathOperator{\im}{Im}

\title[Classification of Lawlor neck]{Stationary integral cycles asymptotic to a pair of transverse oriented $n$-planes at infinity}

\author{Xuanyu Li}
\address{Department of Mathematics, Cornell University, Ithaca, NY 14853, USA}
\email{xl896@cornell.edu}

\author{Zhihan Wang}
\address{Department of Mathematics, Cornell University, 310 Malott Hall, Ithaca, NY 14853}
\email{zw782@cornell.edu}

\date{}

\begin{document}
\pagenumbering{arabic}

\begin{abstract}
  Let $P, Q$ be two oriented transverse $n$-planes in $\R^{2n}$, $n\geq 2$, which satisfy the equality in the Lawlor--Nance angle criterion. We prove that any stationary integral $n$-cycle in $\R^{2n}$ with $[P]+[Q]$ being the tangent cone at infinity, in both the current and varifold senses, is a translation of $[P]+[Q]$ or of a Lawlor neck. 
\end{abstract}

\maketitle

\setcounter{tocdepth}{1}
\tableofcontents

\section{Introduction} \label{Sec_Intro}

  Given a pair of oriented $n$-planes $P, Q$ in $\R^{2n}$, $n\geq 2$, there are angles $0\leq \theta_1\leq \cdots\leq \theta_{n-1}\leq \theta_n\leq \pi-\theta_{n-1}$ uniquely determined by $P, Q$ and some orthonormal basis $e_1,\ldots,e_n,f_1,\ldots, f_n$ of $\R^{2n}$ such that \[
   \xi_P = e_1\wedge\cdots\wedge e_n, \qquad
   \xi_Q = \bigwedge_{j=1}^n(\cos\theta_j e_j+\sin\theta_j f_j),
  \]
  where $\xi_P, \xi_Q$ are the oriented unit $n$-vectors of $P, Q$. 
  Note that $P, Q$ are transverse to each other if and only if $\theta_1>0$.
  
  Nance \cite{Nance87} and Lawlor \cite{Lawlor89} showed that as integral $n$-cycles, $[P]+[Q]$ is area-minimizing if and only if the following angle criterion holds: 
  \begin{align}
    \theta_n \leq \theta_1 + \dots + \theta_{n-1} \,,   \label{eq:angle-ineq} 
  \end{align}

  The equality case is particularly interesting: If 
  \begin{align}
    \theta_n = \theta_1 + \dots + \theta_{n-1} \,,   \label{eq:angle-condition} 
  \end{align}
  then besides the planar cone, $[P]+[Q]$ admits a family of smooth embedded connected area-minimizing submanifold $\cL$, known as the \emph{Lawlor necks} \cite{Lawlor89}, asymptotic to it near infinity (see Section \ref{sec:dimension-two} and Section \ref{sec:model} for its definition when $n=2$ and $n\geq 3$ respectively). 
  Moreover, under some orthogonal almost complex structure on $\R^{2n}$, both $P\cup Q$ and $\cL$ are special Lagrangian with the same phase. 
  
  A natural question is to classify general complete oriented minimal submanifolds $M\subset \R^{2n}$ with tangent cone $[P]+[Q]$ at infinity in both the current and varifold sense. In this direction, Imagi, Joyce, and Oliveira dos Santos \cite{IJS16} proved the uniqueness of such $M$ up to translation and dilation among smooth and exact special Lagrangian submanifolds using Floer theory. Their work is generalized to special Lagrangian integral $n$-cycles under a certain smooth exact approximation condition by Li and Sz\'ekelyhidi \cite{LiSzekelyhidi24}. 
  In a more recent work, Lee \cite{Lee26} proves the classification of such minimal submanifolds $M$ under the assumption $H^1(M)=0$, which allows him to reduce to the classification result by Imagi, Joyce and Oliveira dos Santos \cite{IJS16}; he also proves the classification of any stationary integral $n$-cycles asymptotic to $[P]+[Q]$ when the angles between $P, Q$ are $\theta_1=\dots =\theta_{n-1}=\pi/n$ and $\theta_n=(n-1)\pi/n$.

  The goal of this paper is the following classification in general:
  \begin{theorem}\label{thm:main}
    Let $n\geq 2$, $P, Q\subset\R^{2n}$ be transverse oriented $n$-planes satisfying \eqref{eq:angle-condition}. Let $T$ be a stationary integral $n$-cycle in $\R^{2n}$. Suppose that, with $\eta_R(x):=x/R$,
    \begin{equation}\label{eq:blowdown}
      (\eta_R)_\# T\to[P] + [Q], \qquad
      (\eta_R)_\# |T|\to |P|+|Q| \qquad \text{ when }R\to +\infty
    \end{equation}
    locally as currents and varifolds. Then there exists $c\in\R^{2n}$ such that \[
      \text{ either }\  T= [c + P] + [c + Q], \qquad\text{ or }\ 
      T= [c + \cL] \,;
    \]
    where $\cL$ is a Lawlor neck asymptotic to $[P]+[Q]$ at infinity. 
  \end{theorem}

  The notion of stationary integral cycles are introduced in Section \ref{sec:calibration}. In particular, every $n$-dimensional oriented properly immersed smooth minimal submanifold in $\R^{2n}$ is naturally a stationary integral cycle, and the Theorem applies to it.

  One may naturally ask for classifications when the strict inequality holds in \eqref{eq:angle-ineq}. We shall address this in an upcoming work.


\subsection{Sketch of the proof}
  Let $T$ be as in Theorem \ref{thm:main}. We first prove in Section \ref{sec:calibration}, following Lee's argument \cite[Theorem 3.5]{Lee26}, that under some orthogonal almost complex structure $J$ on $\R^{2n}$, $T$ is special Lagrangian. 
  When $n=2$, this means $T$ is holomorphic under some different orthogonal almost complex structure, and the classification is well-studied.
  
  For the rest of the introduction, let $n\geq 3$. We can choose a complex coordinate on $(\R^{2n}, J)$ such that 
  \begin{align}
    P = 0\oplus i\R^n, \qquad 
    Q = \{i\Phi(x): x\in \R^n\} = \{x+iAx: x\in \R^n\} \,.  \label{Equ_Intro_P,Q}
  \end{align}
  where $\Phi:= \diag(e^{i\phi_1}, \dots, e^{i\phi_n})$, $A:= -\diag(\cot\phi_1, \dots, \cot\phi_n)$.
  (Here we set $\phi_j:=\theta_j$, $1\leq j\leq n-1$, and $\phi_n:= \pi-\theta_n$.) 
  Let $\Omega:=dz^1\wedge\dots\wedge dz^n$ be the complex volume form. Then $P, Q, T$ are all calibrated by $\re(i^{-n}\Omega)$.
  Moreover, after a translation, outside a large ball, we can write $T$ as the gradient graph of some scalar functions $g_P, g_Q$ over $P, Q$ respectively, and they have the following asymptotic expansion: there are $A_P, A_Q\in \R$ such that for $E\in \{P, Q\}$, \[
    g_E(t) = A_E |t|^{2-n} + O(|t|^{1-n}), \qquad \text{ as } t\to \infty \text{ on }E \,.
  \]
  And as is argued in Section \ref{Subsec_positive-asymp}, let $Z$ be the real radial vector field on $\R^{2n}$, then by integrating $d(\iota_Z\im(i^{-n}\Omega)) = n\im(i^{-n}\Omega)$ over $T\llcorner B_R$ and send $R\to +\infty$, we can conclude that $A_P+A_Q=0$.

  One key ingredient of the proof is then a trace inequality (Lemma \ref{Lem_Tr-ineq}) over Lagrangian planes: Let $A$ is a symmetric $n\times n$ real matrix such that \[
    \sum_{j=1}^n \arctan \lambda_j(A) = \frac{-(n-2)\pi}{2}, \qquad 
    M_A:= \begin{bmatrix}
      2A & -I_n \\ -I_n & 0
    \end{bmatrix} \,.
  \]
  Then for every Lagrangian plane $E\subset \R^n\oplus i\R^n$ with phase $\frac{n\pi}{2}+\pi\Z$, we have $\tr_E M_A \leq 0$. More importantly, equality holds if and only if $E=\{x+iAx: x\in \R^n\}$ or $0\oplus i\R^n$. 
  This is the lower-phase counterpart of Yuan's convexity theorem \cite[Lemma 2.1 \& Step 3]{Yuan06}, we give a direct proof in Section \ref{Subsec_Pf-trace-ineq}.

  As a first application, in Section \ref{Subsec_positive-asymp}, we consider the quadratic function $q(x, y):= \langle x, y\rangle - \langle Ax, x \rangle$, where $A$ is the graphical matrix of $Q$ over $\R^n\oplus 0$ as in \eqref{Equ_Intro_P,Q}. The trace inequality yields $\Delta_T q = -\tr_T M_A\geq 0$ on $\Reg T$, while the asymptotic analysis near infinity proves $q(z)\to 0$ as $z\to \infty$ along $\spt T$. Thus $q\leq 0$. And by strong maximum principle, 
  \begin{itemize}
    \item either $q|_{\spt T} < 0$, and then as $R\to \infty$, by the asymptotic analysis, \[
      0 < \int_{B_R} \Delta_T q\ d\|T\| = \int_{\spt T\cap \partial B_R} \langle \nabla_T q, \nu_R\rangle\ d\cH^{n-1} = c_n(A_Q-A_P) + o(1)
    \]
    for some $c_n>0$. In particular this forces $A_P=-A_Q<0$;
    \item or $q\equiv 0$, which forces the tangent planes of $T$ belonging to $\{P, Q\}$ by the equality case of trace inequality, and thus $T$ is a translation of $[P]+[Q]$.
  \end{itemize}
  We shall proceed under the first alternative. The sign $A_P=-A_Q<0$ allow us to rescale $T$ so that the leading coefficients $A_P, A_Q$ for $T$ agrees with those coefficients of a fixed Lawlor neck $\hat\cL$  asymptotic to $P\cup Q$ near infinity. 

  The next key ingredient is the observation (Proposition \ref{prop:model-graph}) that, the Lawlor neck $\hat\cL$ is a global gradient graph over $(\R^n\setminus \{0\})\times \{0\}$ under the coordiantes above, i.e. for some $U\in C^\infty(\R^n\setminus \{0\})$,  \[ 
    \hat\cL = \{x+ i \nabla U(x): x\in \R^n\setminus \{0\}\} \,.
  \]
  Making use of this $U$, in Section \ref{sec:comparison}, we consider a different potential $\hat q(x, y):= \langle x, y\rangle - 2U(x)$. The special Lagrangian equation for $U$ allows us to use the trace inequality again to conclude that $\Delta_T \hat q = -\tr_T M_{\nabla^2 U} \geq 0$ on $\Reg T$, and $\Delta_{\hat\cL}\hat q|_{\hat\cL}\equiv 0$. 
  On the other hand, the matching asymptotics for $T$ and $\hat \cL$ near infinity yields via integration by part \[
    \lim_{R\to +\infty}\int_{B_R} \Delta_T\hat q\ d\|T\| 
    = \lim_{R\to +\infty}\int_{B_R} \Delta_{\hat\cL}\hat q\ d\|\hat\cL\| = 0\,,
  \]
  which implies $\Delta_T\hat q\equiv 0$. In addition, the rigidity case in the trace inequality implies that the tangent plane at every $x+iy\in \Reg T$ is equal to the tangent plane of $\hat\cL$ at $x+i\nabla U(x)$, thus concluding the proof that $T$ is a translation of a Lawlor neck.

\subsection*{Acknowledgment}
  The authors are grateful to Prof. Xin Zhou and Prof. Daniel Stern for their continuous support and valuable guidance. 
  Xuanyu Li is partly supported by Simons Dissertation Fellowship in Mathematics.

\subsection*{AI disclosure statement}
  The second named author started to work on this classification problem at the beginning of 2026, and he had a rough idea to derive the classification using min-max in August 2026. 
  After ChatGPT 6 came out in September, the authors used the problem to test the ability of AI and obtained a proof purely generated by ChatGPT Astra. The authors then identified the two key ingredients in the proof as summarized in the introduction.  
  We think this is a more elegant approach, and the purpose of these notes is therefore to organize and explain this argument independently. 
  No AI-generated text was included.

%

\section{Preliminaries} \label{sec:calibration}
  We recall basic terms and facts from geometric measure theory; interested readers may refer to further background in \cite{Simon83_GMT}. 

  Let $N>n\geq 2$. An \emph{integral $n$-cycle} $T$ in $\R^N$ is an $n$-dimensional integer rectifiable current in $\R^N$ with locally finite mass and $\partial T = 0$.
  For an oriented submanifold $M\subset \R^N$, we write $[M]$ for its associated integral current with multiplicity $1$. 
  For an integral current $T$, let $|T|$ denote its associated integral varifold, $\|T\|$ denote the mass measure of $T$ and $\xi_T$ denote its oriented unit approximate tangent $n$-vector.
  
  We call $T$ \emph{stationary} if the first variation of the varifold $|T|$, $\delta|T|(V)=0$ for every compactly supported smooth ambient vector field $V$ on $\R^{N}$. 

\subsection{Complex coordinates from a pair of plane} \label{Subsec_Cplx-coordinates}
  Throughout the paper, let $n\geq 2$, $P, Q$ be $n$-dimensional oriented planes in $\R^{2n}$ with $P\cap Q = 0$. By \cite{Morgan83}, there are angles $0< \theta_1\leq \cdots\leq \theta_{n-1}\leq \theta_n\leq \pi-\theta_{n-1}$ uniquely determined by $P, Q$ and oriented orthonormal frame $\{e_j\}_j$, $\{f_j\}_j$ of $P, P^\perp$ respectively, such that \[
   \xi_P = e_1\wedge\cdots\wedge e_n, \qquad
   \xi_Q = \bigwedge_{j=1}^n(\cos\theta_j e_j+\sin\theta_j f_j) \,.
  \]
  $\{\theta_j\}$ are called the \emph{characterizing angles} of $P, Q$.
  And let $\phi_j\in (0, \pi)$ be the adapted angles given by $\phi_j:= \theta_j$ for $1\leq j\leq n-1$ and $\phi_n:= \pi-\theta_n$.
  
  Define an orthogonal almost complex structure $J$ on $\R^{2n}$ by \[
    J e_j := f_j, \quad 1\leq j\leq n-1\,; \qquad
    J e_n := -f_n \,.
  \]
  Also choose a complex coordinate of $(\R^{2n}, J)$ by setting the real coordinate vectors to be $-Je_j$ and the imaginary coordinate vectors to be $e_j$. 
  Note that under this choice, we identify $\R^{2n}=\R^n\oplus i\R^n = \C^n$, and \[
    P = 0\oplus i\R^n, \qquad 
    Q = \{i\Phi(x): x\in \R^n\}, \  \text{ where }\Phi:= \diag(e^{i\phi_1}, \dots, e^{i\phi_n}) \,.
  \] 
  as unoriented $n$-planes. 
  We emphasize that the orientation of $P$ agrees with the one induced by $\R^n\mapsto P\subset \C^n,\ x\mapsto ix$, but the orientation of $Q$ is the negative of the induced orientation by $i\Phi: \R^n\to Q\subset \C^n$, i.e. $(i\Phi)^*\xi_Q = -dx^1\wedge \dots \wedge dx^n$.  

  Throughout the paper, we also let \[
    \Omega:= dz^1\wedge \dots \wedge dz^n\,, \qquad  
    \omega:= dx^1 \wedge dy^1 + \dots + dx^n \wedge dy^n 
  \]
  be the holomorphic volume form and symplectic form of $(\R^{2n}, J)$.
  
  Suppose now $P, Q$ are from Theorem \ref{thm:main}, and hence $\phi_1+\dots + \phi_n = \pi$. By Harvey-Lawson \cite{HL82}, with the above convention of orientation, $P, Q$ are both Lagrangian planes (i.e. $\omega|_P=\omega|_Q = 0$) and both calibrated by the constant real $n$-form $\re(i^{-n}\Omega)$, i.e. \[
    \langle \re(i^{-n}\Omega), \xi_P \rangle = \langle \re(i^{-n}\Omega), \xi_Q \rangle = 1\,; \qquad 
    |\langle \re(i^{-n}\Omega), \xi_{L} \rangle|\leq 1 \ \ \forall\, \text{oriented $n$-plane }L\subset \R^{2n} \,.
  \] 
  
  In general, call an oriented Lagrangian plane $L\subset \C^n$ have \emph{(oriented) phase $\al\in \R/2\pi\Z$} if $\langle \re(e^{-i\al}\Omega), \xi_L \rangle = 1$ (which also implies that $\im(e^{-i\al}\Omega)|_L = 0$). Under the convention above, both $P, Q$ have phase $n\pi/2$.
  When $L$ is unoriented, any choice of orientation well-defines the phase modulo $\pi$.

\subsection{Calibrated cycle}
  We now adapt the argument of Lee \cite[Section 3]{Lee26} to show that $T$ is calibrated by $\re(i^{-n}\Omega)$ (following Harvey-Lawson \cite{HL82}, we call such $T$ \emph{special Lagrangian}).
  
  \begin{lemma} \label{lem:structure}
    Let $P, Q, T$ be as in Theorem~\ref{thm:main}, and endow $\R^{2n}$ with the complex structure and coordinates in Section \ref{Subsec_Cplx-coordinates}. Then $T$ is calibrated by $\re(i^{-n}\Omega)$.

    Moreover, when $n\geq 3$, there's a unique $c\in \R^{2n}$ (called the ``center" of $T$) such that outside a large ball, the translated cycle $T-c$ is the disjoint union of normal gradient graphs over $P, Q$, i.e. 
    for $E\in \{P, Q\}$, there exists $g_E\in C^\infty(E)$ such that for $R\gg 1$, \[
      Z_P(P)\cup Z_Q(Q)\setminus B_R = \spt (T-c)\setminus B_R\,, \qquad
      \text{ where } Z_E(t) := t + J\nabla_E g_E(t), \ \ t\in E \,;
    \]
    and $g_E$ has an asymptotic expansion
    \begin{align}
      g_E(t) = A_E|t|^{2-n}+O_\infty(|t|^{1-n}) \qquad \text{ as }t\to \infty \text{ on }E \,.   \label{Equ_Prelim_g_E-expansion}
    \end{align}
    for some constant $A_E\in\R$. Here and below, $O_\infty(|t|^\lambda)$ stands for a function whose $k$-th order derivative decays like $O(|t|^{\lambda-k})$ for every $k\geq 0$.
  \end{lemma}
  \begin{proof}
    Since the tangent cone at infinity of $T$ as a varifold is a multiplicity one regular cone $|P|+|Q|$, with integrable link (disjoint union of equatorial $S^{n-1}$), by the asymptotic analysis of Allard-Almgren \cite{AA81}, outside a compact subset, $\spt T$ is the disjoint union of normal graphs of some $u_E: E\setminus B_{R_0}\to E^\perp$, $E\in\{P, Q\}$, with expansion
    \begin{equation}\label{eq:minimal-expansion}
      u_E(t) =  \begin{cases}
        a_E + b_E \log |t| + O_\infty(|t|^{-1}), & n = 2,\\
        a_E + b_E |t|^{2-n} + O_\infty(|t|^{1-n}), & n\geq 3,
      \end{cases}
    \end{equation}
    where $a_E, b_E\in E^\perp$. Following the term \cite[Definition 3.2]{Lee26}, this means $T$ is regular at infinity. Since $P, Q$ span $\R^{2n}$, there's a unique $c\in \R^{2n}$ such that $a_P=a_Q=0$ in the expansion associated to the translated cycle $T-c$. 
    Also, by \cite[Theorem 3.1]{Lee26}, $\Sing T$ has Hausdorff dimension $\leq n-2$.
    Then by \cite[Theorem 3.5]{Lee26}, since $P, Q$ are both calibrated by $\re(i^{-n}\Omega)$, so is $T-c$ (and thus $T$), and $b_P=b_Q=0$ in the expansion associated to $T-c$. In particular, $\Reg T$ is a Lagrangian submanifold.
    
    When $n\geq 3$, since the graph of $u_E$ over $\R^n\setminus B_{R_0}$ is Lagrangian and simply connected, we have $u_E = J\nabla g_E$ for some $g_E\in C^\infty(E\setminus B_{R_0})$ (extended arbitrarily to the whole $E$). 
    The special Lagrangian equation for $g_E$ (from the fact $Z_E^*\im(i^{-n}\Omega) = 0$) has linear term $\Delta g_E$ and error terms of degree $\geq 2$. Thus the asymptotic \eqref{eq:minimal-expansion} gives the expansion of $g_E$ in \eqref{Equ_Prelim_g_E-expansion} near infinity.
  \end{proof}

  \begin{remark} \label{Rem_T_R,conormal}
    We record some useful consequences of the asymptotic expansion \eqref{Equ_Prelim_g_E-expansion}, for $n\geq 3$ and with $T$ translated to its center as in Lemma~\ref{lem:structure}.
    First, the pull-back metrics by $Z_E$ from $\R^{2n}$ satisfy
    \begin{equation}\label{eq:end-metrics}
      G_E := Z_E^*\ g_\Euc = I_n + (D^2g_E)^2 = I_n + O_\infty(|t|^{-2n}),\qquad E \in \{P, Q\}\,.
    \end{equation}
    For $R\gg 1$, let $T_R$ be the compact truncation of $T$ obtained by cutting each end at the image of $E\setminus B_R$ under parametrization $Z_E$. In this way, $\spt\partial T_R = \partial_P T_R\cup \partial_Q T_R$, where \[
      \partial_E T_R := Z_E(\{t\in E: |t|=R\}),  \qquad E\in\{P,Q\} \,,
    \]
    and let $\nu_R$ denote the outward unit co-normal along these boundary components in $\Reg T$. Then $\nu_R$ and the area measure of $\partial_E T_R$ are given by 
    \begin{equation}\label{eq:conormal-density}
      \begin{aligned}
        \nu_R\circ Z_E(t)
          & = \langle G_E^{-1}(t)t,t\rangle^{-1/2}\cdot DZ_E(t)[G_E^{-1}(t)t] \,,  \\
        Z_E^*\bigl(d\cH^{n-1}|_{\partial_E T_R}\bigr)
          & = R^{-1}\sqrt{\det G_E(t)}\, \sqrt{\langle G_E^{-1}(t)t,t\rangle}\ d\cH^{n-1}|_{\partial B_R}\,,
      \end{aligned}
    \end{equation}
    Indeed, the first identity is the normalized push-forward of the $G_E$-gradient of $|t|$, and the second follows from the coarea formula for $|t|$.
  \end{remark}

  We have seen in the proof that $\Sing T$ has Hausdorff dimension $\leq n-2$.
  The calibration $\re(i^{-n}\Omega)$ makes $T$ locally area minimizing \cite{HL82}. 
  Also, since $T$ has a multiplicity $1$ tangent cone with an isolated singularity at infinity, Allard regularity \cite{Allard72} guarantees the compactness of $\Sing T$. 
  We have the following dichotomy about $\Reg T$.
  \begin{corollary} \label{prop:components}
    Let $P, Q, T$ be in Theorem \ref{thm:main}.
    Then after possibly a translation of $T$, we have, either $T=[P]+[Q]$, or $\Reg T$ is connected.
  \end{corollary}
  \begin{proof}
    Let $M_1, \dots, M_k$ be the connected components of $\Reg T$, each is an oriented $n$-submanifold in $\R^{2n}$ calibrated by $\re(i^{-n}\Omega)$. Moreover, $\overline{M_j}\setminus M_j \subset \Sing(T)$, which has dimension $\leq n-2$, and $\cH^n\llcorner M_j$ is bounded by $\|T\|$ and hence locally finite in $\R^{2n}$. Therefore, the associated currents $[M_j]$ are minimizing $n$-cycles in $\R^{2n}$, and $\sum_j\|M_j\|\leq \|T\|$. By area-monotonicity formula \cite{Simon83_GMT}, \[
      \lim_{R\to +\infty} (\omega_n R^n)^{-1} \|M_j\|(B_R) \geq 1 \,.
    \]
    Hence the asymptotic assumption on $T$ forces $k\leq 2$; and $k=2$ implies \[
      \lim_{R\to +\infty} (\omega_n R^n)^{-1} \|M_j\|(B_R) = 1, \quad j=1,2 \,; 
    \]  
    which implies that $M_j$ are both planes by Allard regularity \cite{Allard72}, and forces $T=[P]+[Q]$ after a translation. 
    And $k=1$ means $\Reg T$ is connected.
  \end{proof}

\subsection{Classification in dimension $n=2$}\label{sec:dimension-two}
  In dimension $n=2$, a special Lagrangian cycle in $\C^2$ is a holomorphic curve under a different orthogonal almost complex structure.  We thus have the following direct classification result.
  \begin{proposition}\label{prop:dimension-two}
    Let $n=2$ and $P, Q, T$ be in Theorem \ref{thm:main}.
    Then there are an orthogonal complex structure $I$ and linear independent complex-linear functions $\ell_P,\ell_Q:(\R^4, I)\to \C$ with kernels $P,Q$, such that, after a translation, \[
      T=[P]+[Q] \qquad \text{or} \qquad
      T=[\{\ell_P\ell_Q=c\}] \ \ \text{ for some } c\in\C\setminus\{0\} \,.
    \]
    $\{\ell_P\ell_Q=c\}$ in the second alternative are called \emph{Lawlor necks}.
  \end{proposition}
  \begin{proof}
    This is a consequence of \cite[Theorem 6.8]{LambertLotaySchulze21}.
  \end{proof}

\section{Trace inequality and the leading order asymptotic} \label{Sec:Tr-ineq}
  The central Lemma of this section is the following trace inequality. Throughout the section, let $n\geq 3$, $A$ be a real symmetric $n\times n$ matrix. Define 
  \begin{align}
    F(A) := \sum_{j=1}^n \arctan\lambda_j(A)\in (-\frac{n\pi}2, \frac{n\pi}{2})\,, & &
    M_A := \begin{bmatrix}
      2A & -I_n \\ -I_n & 0
    \end{bmatrix} \,.  \label{Equ_Def_F(A),M_A}
  \end{align}
  Note that $F(A)+\pi\Z$ is the phase of the graph $G_A:=\{x+iAx:x\in \R^n\}\subset \C^n$  (here we only work with phase for unoriented Lagrangian plane, which is well-defined$\mod \pi$).
  \begin{lemma}[Trace inequality]\label{Lem_Tr-ineq}\label{cor:graph-vertical}
    Let $\eta\in \R$, with phase understood modulo $\pi$. Then the following are equivalent:
    \begin{enumerate}[label={\normalfont(\roman*)}]
      \item\label{Item_Tr-ineq_tr_E M_A<0} $\tr_E M_A\leq 0$ for every Lagrangian plane $E\subset \R^{2n}=\R^n\oplus i\R^n$ with phase $\eta+\pi\Z$.
      \item\label{Item_Tr-ineq_phase E} $(\eta+\pi\Z)\cap (-\frac{n\pi}2, F(A))= \emptyset$.
    \end{enumerate}
    Moreover, for $\eta$ in the range of \ref{Item_Tr-ineq_phase E}, equality holds for $E$ in \ref{Item_Tr-ineq_tr_E M_A<0} if and only if $E=G_A$ (with phase $F(A)+\pi\Z$) or $0\oplus i\R^n$ (with phase $\frac{n\pi}{2}+\pi\Z$). 
  \end{lemma}

  For a graphical Lagrangian plane $E=G_B$ in the lower branch $F(B)\leq -(n-2)\pi/2$, the inequality in \ref{Item_Tr-ineq_tr_E M_A<0} is the supporting-hyperplane inequality for the convex sublevel set $\{F\leq F(B)\}$. This is the lower-phase counterpart of Yuan's convexity theorem \cite[Lemma 2.1 \& Step 3]{Yuan06}: condition \ref{Item_Tr-ineq_phase E} gives $F(A)\leq F(B)$, and the matrix inequality follows by the doubly stochastic argument below. We give a direct proof in Section \ref{Subsec_Pf-trace-ineq}, and carefully study the equality case. 
  In Section \ref{Subsec_positive-asymp}, we give a first application of it, which proves that under the setup of Lemma \ref{lem:structure}, the leading coefficients $A_P, A_Q$ as in \eqref{Equ_Prelim_g_E-expansion} satisfy $A_P=-A_Q<0$ unless $T$ is a translation of $[P]+[Q]$.

\subsection{Proof of Lemma \ref{Lem_Tr-ineq}.} \label{Subsec_Pf-trace-ineq} 
  We first compute the trace over each graph \[
    G_B:=\{x+iBx: x\in \R^n\}\subset \R^n\oplus i\R^n
  \] 
  where $B$ is an arbitrary real $n\times n$ symmetric matrix:  Take $\{e_k\}_k$ an orthonormal basis of $\R^n$ such that $Be_k = b_ke_k$, where $b_k\in \R$ are eigenvalues of $B$, then $\{(1+b_k^2)^{-1/2}(e_k+ib_ke_k)\}_k$ is an orthonormal basis of $G_B$, hence 
  \begin{align}
   \begin{split}
    \tr_{G_B} M_A & = \sum_{k=1}^n (1+b_k^2)^{-1}\cdot \begin{bmatrix} e_k^T & b_ke_k^T \end{bmatrix} \begin{bmatrix}
      2A & -I_n \\ -I_n & 0
    \end{bmatrix} \begin{bmatrix} e_k \\ b_ke_k \end{bmatrix} \\
    & = 2\sum_{k=1}^n  (1+b_k^2)^{-1}(\langle A e_k, e_k\rangle - b_k)
    = 2\tr\left((I_n+B^2)^{-1}(A-B) \right) \,.
   \end{split} \label{Equ_Tr-ineq_tr_(G_B)M_A}
  \end{align}
  We denote for simplicity $a_k:=\lambda_k(A)$ to be the eigenvalues of $A$. Taking $B:=A-tI_n$ yields \[
    \tr_{G_{A-tI_n}} M_A = 2\sum_{k=1}^n \frac{t}{1+(a_k-t)^2} > 0, \quad \text{ if }t>0\,.
  \] 
  As $t$ increases from $0$ to infinity, $F(A-tI_n)$ decreases continuously from $F(A)$ to $-n\pi/2$. Thus as long as $(\eta+\pi\Z)\cap (-\frac{n\pi}{2}, F(A))\neq \emptyset$, we can find some $t_\eta>0$ such that $G_{A-t_\eta I_n}$ has phase $F(A-t_\eta I_n) + \pi\Z = \eta+\pi\Z$ but the trace inequality fails. This proves \ref{Item_Tr-ineq_tr_E M_A<0} implies \ref{Item_Tr-ineq_phase E}.

  To prove that \ref{Item_Tr-ineq_phase E} implies \ref{Item_Tr-ineq_tr_E M_A<0}, by the denseness of all graphical Lagrangian planes of a fixed phase in the space of all Lagrangian planes of a fixed phase, it suffices to prove that $\tr_{G_B}M_A\leq 0$ for any symmetric matrix $B$ with $\eta:=F(B)$ satisfying \ref{Item_Tr-ineq_phase E}. 
  Let $A_{kk}:= \langle A e_k, e_k\rangle$, and $a_1, \dots, a_n$ be the eigenvalues of $A$. The following fact allows us to only focus on the case when $A, B$ are simultaneously diagonalized:
  \begin{Claim} \label{Claim_Permutation}
    $(A_{11}, A_{22}, \dots, A_{nn})$ is a convex combination of the permutations of eigenvalues: \[
      \{(a_{\sigma(1)}, a_{\sigma(2)}, \dots, a_{\sigma(n)}): \sigma\in \mfk S_n\} \,.
    \] 
  \end{Claim}
  \begin{proof}
    Let $\Omega\in SO(n)$ such that $A=\Omega^T\diag(a_1, \dots, a_n)\Omega$. Then $A_{kk} = \sum_j\Omega_{jk}^2a_j$, $1\leq k\leq n$. Since $[\Omega_{jk}^2]$ is a doubly stochastic matrix (i.e. it has nonnegative entries, and the sum of each row or column is $1$), by Birkhoff-von Neumann Theorem \cite[Theorem II.2.3, p.~37]{Bhatia97_Matrix}, $[\Omega_{jk}^2]$ is a convex combination of permutation matrices. This proves the Claim.
  \end{proof}

  With this Claim, by \eqref{Equ_Tr-ineq_tr_(G_B)M_A}, it suffices to deal with the case when 
  \begin{align}
    \tr_{G_B} M_A = 2\sum_{k=1}^n  \frac{a_k - b_k}{1+b_k^2} 
    = 2\sum_{k=1}^n \frac{-\cot\al_k + \cot\beta_k}{1+\cot^2\beta_k} 
    = 2\sum_{k=1}^n \frac{\sin\beta_k\cdot\sin(\al_k-\beta_k)}{\sin\al_k}  \,. \label{Equ_Tr-ineq_Tr_B(M_A)=sum h(al,beta)}
  \end{align}
  here we let $\al_k:= \pi/2 + \arctan a_k$, $\beta_k:= \pi/2 + \arctan b_k\in (0, \pi)$. Note that $F(A)=\sum_k\al_k-\frac{n\pi}{2}$, and the phase condition on $\eta+\pi\Z=F(B)+\pi\Z$ in \ref{Item_Tr-ineq_phase E} is equivalent to 
  \begin{align}
    (0, \al_1+\dots + \al_n) \cap (\beta_1+\dots+\beta_n + \pi\Z) = \emptyset \,. \label{Equ_Tr-ineq_phase-restrict-sum al-beta}
  \end{align}
  This in particular forces $\sum_k\al_k\leq \pi$.
  \begin{Claim} \label{Claim_h(al,beta)}
    Let $h(\al, \beta):= \frac{\sin\beta\cdot\sin(\al-\beta)}{\sin\al}$. Then for $\al, \al'>0$ such that $\al+\al'<\pi$, and any $\beta, \beta'\in \R$, there holds 
    \begin{align*}
      h(\al,\beta) + h(\al', \beta') 
      & = h(\al+\al', \beta+\beta') - \frac{|\sin\al' \sin\beta - e^{i(\beta+\beta' - \al-\al')}\sin\al\,\sin\beta'|^2}{\sin\al\,\sin\al'\,\sin(\al+\al')} \\
      & \leq h(\al+\al', \beta+\beta')  \,.
    \end{align*}
    In particular, if $\beta, \beta'\in (0, \pi)$, then the equality holds in the last inequality above if and only if $\beta=\al$ and $\beta'=\al'$.
  \end{Claim}
  \begin{proof}
    The first equality follows from direct expansion, and the denominator in the remainder is positive. 
    If $\beta,\beta'\in(0,\pi)$, vanishing of the remainder forces $e^{i(\beta+\beta'-\al-\al')}=1$. Since $-\pi<\beta+\beta'-\al-\al'<2\pi$, this gives $\beta+\beta'=\al+\al'$. Thus equality holds exactly when \[
      \gamma:=\beta+\beta' = \al+\al'\in (0, \pi), \qquad \sin\al'\sin\beta = \sin\al\, \sin\beta' \,.
    \]
    Notice that this implies \[
      0 = \sin(\gamma-\al)\sin\beta - \sin\al\, \sin(\gamma-\beta) = \sin\gamma\, \sin(\beta-\al)
    \]
    and hence forces $\beta=\al$, $\beta'=\al'$.
  \end{proof}

  Combine this Claim with \eqref{Equ_Tr-ineq_Tr_B(M_A)=sum h(al,beta)} inductively yields (since $\al_k>0$ and $\sum_k\al_k\leq \pi$), \[
    \tr_{G_B} M_A = 2\sum_{k=1}^n h(\al_k, \beta_k) \leq 2h(\al_1+\dots+\al_{n-1}, \beta_1+\dots+\beta_{n-1}) + 2h(\al_n, \beta_n) \,.
  \]
  If $\sum_k\al_k<\pi$, then $\tr_{G_B}M_A \leq 2h(\sum_k \al_k, \sum_k\beta_k)\leq 0$ by \eqref{Equ_Tr-ineq_phase-restrict-sum al-beta}, with $\tr_{G_B}M_A=0$ if and only if $\al_k=\beta_k$, $\forall\, 1\leq k\leq n$. 
  While if $\sum_k\al_k = \pi$, then \eqref{Equ_Tr-ineq_phase-restrict-sum al-beta} forces $\sum_k\beta_k=\ell\pi\in \pi\Z$ and \[
    h(\al_1+\dots+\al_{n-1}, \beta_1+\dots+\beta_{n-1}) = h(\pi-\al_n, \ell\pi-\beta_n) = -h(\al_n, \beta_n)\,,
  \]
  which again proves $\tr_{G_B}M_A\leq 0$. And $\tr_{G_B}M_A=0$ implies $\al_k=\beta_k$, $\forall\, 1\leq k\leq n-1$, which then forces $\al_n=\beta_n$ since $\al_k, \beta_k\in (0, \pi)$.

  For matrices $B$ not simultaneously diagonalized with $A$, equality in the convex combination from Claim \ref{Claim_Permutation} requires every permutation with positive weight to pair equal eigenvalues with $B$. Hence $A_{kk}=b_k$ and $\sum_k a_k^2=\sum_k b_k^2$, so
  \[
    \sum_{j\ne k}|\langle Ae_j,e_k\rangle|^2=\tr(A^2)-\sum_k b_k^2=0.
  \]
  Thus $A=B$ in this case as well.

  In summary, we have shown that for every symmetric $n\times n$ real matrix $B$ whose graph $G_B$ has phase $\eta$ satisfying \ref{Item_Tr-ineq_phase E}, there holds $\tr_{G_B}M_A\leq 0$, with equality if and only if $B=A$. 
\medskip

  Now let $E$ be a general fixed non-graphical Lagrangian plane with phase $\eta$ satisfying \ref{Item_Tr-ineq_phase E}. We approximate $E$ by graphs with the same phase: let $V\subset \R^n$ such that $iV:= E\cap i\R^n$, then $k:=\dim V>0$. The Lagrangian condition gives  $E= iV\oplus \{x+iBx:x\in V^\perp\}$ for a symmetric map $B:V^\perp\to V^\perp$.
  If $k< n$, rotate $V\oplus iV$ by $e^{i\eps}$ and $V^\perp\oplus iV^\perp$ by $e^{-ik\eps/(n-k)}$. These unitary maps have determinant $1$, preserve the phase, and give graphical planes for $0<\eps\ll 1$. 
  If $k=n$, rotate the coordinate axes by $\eps,\dots, \eps, -(n-1)\eps$ instead. 
  In either case, $E$ is a limit of Lagrangian graphical plane of phase $\eta$. By continuity of the trace, we still have $\tr_E M_A\leq 0$. This finishes the proof that \ref{Item_Tr-ineq_phase E} implies \ref{Item_Tr-ineq_tr_E M_A<0}.

  Moreover, if $\tr_E M_A=0$. Then by the trace inequality proved above $A'\mapsto \tr_E M_{A'}$ has a local maximum at $A$ in the hypersurface $\{A':F(A') = F(A)\}$ of the space of $n\times n$ symmetric matrices. Since $DF|_A[I_n]=\tr((I_n+A^2)^{-1})>0$, the level set is smooth. By the Lagrange multiplier theory, there exists $c\in \R$ such that for every symmetric $n\times n$ real matrix $K$, there holds \[
    \left.\frac{d}{dt}\right|_{t=0} \tr_E M_{A+tK} = c DF|_A[K] = c\tr\left((I_n+A^2)^{-1}K \right)\,.
  \]
  Let $W:=B+iC$ be an $n\times n$ matrix whose columns form an orthonormal basis $\{u_k+iv_k\}_{1\leq k\leq n}$ of $E$, where $B, C$ are real matrices and $u_k, v_k$ are real $n$-vectors. Since $E$ is non-graphical, we know that $B$ is not full rank. Then \[
    \left.\frac{d}{dt}\right|_{t=0} \tr_E M_{A+tK} = \sum_{k=1}^n \begin{bmatrix} u_k^T & v_k^T \end{bmatrix} \begin{bmatrix}
      2K & 0 \\ 0 & 0
    \end{bmatrix} \begin{bmatrix} u_k \\ v_k \end{bmatrix}
    = 2\tr(B^TKB) = 2\tr(BB^TK) \,.
  \]
  Combining the two equations above proves that $2BB^T=c(I_n+A^2)^{-1}$. Since $B$ is not full rank, we must have $c=0$, and then $B=0$. Thus $E=0\oplus i\R^n$. 
  This finishes the study of the equality case for a general Lagrangian plane.

\subsection{Leading order asymptotic constraint} \label{Subsec_positive-asymp}
  We give a first application of the trace inequality.
  \begin{proposition} \label{prop:coefficient}
    Let $P, Q, T$ be as in the main Theorem \ref{thm:main}, and suppose $T$ has center $0$ as in Lemma \ref{lem:structure}; let $A_P, A_Q$ be the leading-order coefficients as in \eqref{Equ_Prelim_g_E-expansion}. 

    Then either $T=[P]+[Q]$, or $\spt T\cap \{x=0\}=\emptyset$ and $A_P=-A_Q<0$.
  \end{proposition}
  \begin{proof}
    Under the complex coordinates of Section \ref{sec:calibration}, we have
    \[
      P=0\oplus i\R^n\,, \qquad Q=\{x+iAx:x\in\R^n\}\,;
    \]
    where $A=-\diag(\cot\phi_1,\dots,\cot\phi_n)$ and $\sum_j\phi_j=\pi$. Thus $F(A)=-(n-2)\pi/2$. Both oriented planes are calibrated by $\re(i^{-n}\Omega)$, with the orientations in Section \ref{Subsec_Cplx-coordinates}, and have unoriented phase $n\pi/2+\pi\Z$. 

    We first prove $A_P+A_Q=0$. Let $Z$ be the real radial vector field, $\omega$ the K\"ahler form of $J$, and $\beta:=\iota_Z\im(i^{-n}\Omega)$. Note that
    \begin{align}
      d\beta=\cL_Z\im(i^{-n}\Omega)=n\im(i^{-n}\Omega)\,;\qquad
      \beta|_T=*\bigl(-\iota_Z\omega\bigr)|_T\,.
      \label{Equ_Tr-ineq_dbeta=nIm(Omega)}
    \end{align}
    In fact, the first equation follows from Cartan's formula. To prove the second, take an oriented orthonormal basis $\{e_j\}$ of a calibrated tangent plane. Then
    $i^{-n}\Omega=(e_1^*+i(Je_1)^*)\wedge\dots\wedge(e_n^*+i(Je_n)^*)$.
    Writing $Z^\perp=\sum_j b_jJe_j$, we have
    \[
      \beta|_T=\im\bigg(\sum_j(-1)^{j-1}b_j i\,
        e_1^*\wedge\dots\wedge\widehat{e_j^*}\wedge\dots\wedge e_n^*\bigg)
      =*\bigg(\sum_j b_je_j^*\bigg)
      =*\bigl(-\iota_Z\omega\bigr)|_T\,.
    \]
    This proves the second equality in \eqref{Equ_Tr-ineq_dbeta=nIm(Omega)}.

    For $R\gg1$, let $T_R$ retain the compact core and truncate both ends at base radius $R$, as in Remark \ref{Rem_T_R,conormal}. By Stokes' theorem,
    \[
      0=\int_{T_R}n\im(i^{-n}\Omega)=\int_{\partial T_R}\beta
        =\int_{\partial T_R}*\bigl(-\iota_Z\omega\bigr)|_T\,.
    \]
    Outside a large ball $T$ is the union of the gradient graphs of $g_P,g_Q$. Write $Z_E(t):=t+J\nabla_Eg_E(t)$ intrinsically for the end over $E\in\{P,Q\}$. Then \eqref{Equ_Prelim_g_E-expansion} gives
    \[
      Z_E^*(\iota_Z\omega)
      =d\bigl(\langle t,\nabla_Eg_E(t)\rangle-2g_E(t)\bigr)
      =(n-2)nA_E|t|^{1-n}d|t|+O(|t|^{-n})\,.
    \]
    Let $\partial_ET_R$ denote the boundary component on end $E$. Using the metric estimate \eqref{eq:end-metrics}, we obtain
    \begin{align*}
      0=\int_{\partial T_R}*\bigl(-\iota_Z\omega\bigr)|_T
       &=-\sum_{E\in\{P,Q\}}\int_{\partial_ET_R}(n-2)nA_E R^{1-n}*dr+o(1)\\
       &=(2-n)n|\Sph^{n-1}|(A_P+A_Q)+o(1)\,.
    \end{align*}
    Here $r$ is the base radial function on each end. The reversed orientation on the $Q$ end reverses both the Hodge star and the boundary orientation, so the signs cancel in the scalar outward flux. Sending $R\to+\infty$ proves $A_P+A_Q=0$.
\medskip

    Now suppose $T\ne[P]+[Q]$. Then $\Reg T$ is connected by Corollary~\ref{prop:components}. Consider
    \begin{equation}\label{eq:q}
      q(x,y):=\langle x,y\rangle-\langle Ax,x\rangle.
    \end{equation}
    This is the quadratic function from the proof sketch, expressed in the present coordinates. Since $\nabla^2q=-M_A$, Lemma~\ref{Lem_Tr-ineq} gives on $\Reg T$
    \begin{equation}\label{eq:q-trace}
      \Delta_Tq=-\tr_T M_A\geq0\,;\qquad
      \Delta_Tq(z)=0\ \Longleftrightarrow\ \Tan_zT\in\{P,Q\}.
    \end{equation}
    Also $q|_P=q|_Q=0$, and \eqref{Equ_Prelim_g_E-expansion} gives $\dist(z,P\cup Q)=O(|z|^{1-n})$ along the ends. Hence $q=O(|z|^{2-n})= o(1)$ at infinity. By the weak maximum principle, $q\leq 0$ on $\Reg T$. (The weak maximum principle applies to $T$ even when $\Sing T \neq \emptyset$, since one can consider integrating $(q-\eps)_+\cdot \Delta_T q$ over $T$ and conclude $|\nabla(q-\eps)_+|\equiv 0$ for every $\eps>0$.)
    
    If $q=0$ at some point in $\Reg T$, the strong maximum principle then gives $q\equiv 0$. Equation~\eqref{eq:q-trace} forces every regular tangent plane to be $P$ or $Q$, and by continuity only one choice occurs, contradicting the two-ends assumption. Therefore $q<0$ on $\Reg T$. Since $T$ is smooth near infinity, for $R\gg 1$, we have $\sup_{\partial B_R \cap \spt T} q< 0$, applying weak maximum principle as above again proves $q<0$ on $\spt T$, and hence $\spt T\cap \{x=0\} = \emptyset$.
    
    The same equality characterization shows that $-\tr_T M_A$ is not identically zero on $\Reg T$. Thus, by the first variation formula,
    \[
      0 < \lim_{R\to+\infty}\int_{T_R}(-\tr_T M_A)\ d\|T\|
        = \lim_{R\to+\infty}\int_{\partial T_R}\langle\nabla q,\nu_R\rangle\ d\cH^{n-1}\,;
    \]
    where $\nu_R$ is the outward unit conormal. To calculate this limit, again parametrize the end of $T$ by $Z_E(t)=t+J\nabla_Eg_E(t)$, $E\in \{P, Q\}$. \eqref{eq:conormal-density} and \eqref{eq:end-metrics} give, for $|t|=R$,
    \[
      \nu_R\circ Z_E(t)
       =R^{-1}\bigl(t+J\nabla_E^2g_E(t)[t]\bigr)+O(R^{-2n})\,.
    \]
    Since $\nabla q(z)=-M_Az$, we have
    \begin{align*}
      \langle\nabla q,\nu_R\rangle|_{Z_E(t)}
       &=\left\langle -M_A\bigl(t+J\nabla_Eg_E(t)\bigr),
         R^{-1}\bigl(t+J\nabla_E^2g_E(t)[t]\bigr)\right\rangle
         +O(R^{1-2n})\\
       &=R^{-1}\left\langle JM_At,
           (2-n)^2A_E|t|^{1-n}\nabla_E|t|\right\rangle+O(R^{-n})\,.
    \end{align*}
    Here we use $\langle M_At,t\rangle=0$ since $t\in P\cup Q$, and \eqref{Equ_Prelim_g_E-expansion} gives
    \begin{align*}
      \nabla_Eg_E(t)&=(2-n)A_E|t|^{1-n}\nabla_E|t|+O(|t|^{-n})\,,\\
      \nabla_E^2g_E(t)[t]&=(2-n)(1-n)A_E|t|^{1-n}\nabla_E|t|+O(|t|^{-n})\,.
    \end{align*}
    For $t=(0,y)\in P$ we have $JM_At=-t$, whereas for $t=(x,Ax)\in Q$ we have $JM_At=t$. Consequently
    \begin{align*}
      0&<\lim_{R\to+\infty}\left(
         -\int_{\partial_PT_R}(2-n)^2A_P R^{1-n}\,d\cH^{n-1}
         +\int_{\partial_QT_R}(2-n)^2A_Q R^{1-n}\,d\cH^{n-1}\right)\\
       &=(2-n)^2|\Sph^{n-1}|(A_Q-A_P).
    \end{align*}
    Combined with $A_P+A_Q=0$, this proves $A_P=-A_Q<0$.
  \end{proof}

\section{The Lawlor neck as a global gradient graph}\label{sec:model}

  Let $n\geq 3$. We first recall the classical definition of a Lawlor neck.
  Let $\phi=(\phi_1, \dots, \phi_n)\in (0, +\infty)^n$ such that $\phi_1+\dots+\phi_n = \pi$, and fix $\Lambda>0$.
  In Lawlor's parametrization \cite{Lawlor89} (see also the exposition in \cite[Example~3.9]{IJS16}), a Lawlor neck $\cL_{\phi,\Lambda}\subset \C^n$ is the image of the map
  \begin{equation}\label{eq:lawlor-parametrization}
    F_{\phi, \Lambda}: \R\times \Sph^{n-1} \to \C^n, \quad (s, u) \mapsto (z_1, \dots, z_n) := (\gamma_1(s)u_1, \dots, \gamma_n(s)u_n) \,;
  \end{equation}
  where $\{\gamma_j(s)\in \C\}_{j=1}^n$ is determined by the equations:
  \begin{align}
    F(s) & :=\prod_{j=1}^n(1+s^2/\rho_j), &
    \cP(s) & := (F(s)-1)/s^2,   \label{eq:lawlor-polynomial}\\
    \psi_j(s) & := \int_{-\infty}^s\frac{dt}{(\rho_j+t^2)\sqrt{\cP(t)}}, &
    \gamma_j(s) & :=\sqrt{\rho_j+s^2}\,e^{i\psi_j(s)} \,; \label{eq:lawlor-functions}
  \end{align}
  where the parameters $\{\rho_j>0\}_j$ are the unique positive real numbers such that
  \[
    \psi_j(+\infty)=\phi_j, \quad 1\leq j\leq n\,; \qquad
    \int_{-\infty}^{+\infty} \frac{dt}{2\sqrt{\cP(t)}} = \Lambda \,.
  \]
  We use the convention $\cP(0):=\sum_j\rho_j^{-1}$ to make $\cP$ continuous.

  Orient $\cL_{\phi, \Lambda}$ by setting $d\vol_{\cL_{\phi, \Lambda}}:=\re\Omega|_{\cL_{\phi, \Lambda}}$.
  Note that if $\lambda>0$ is a scaling factor, then
  \[
    F_{\phi, \lambda^2\Lambda}(s, u) = \lambda F_{\phi, \Lambda}(\lambda^{-1}s, u), \qquad
    \cL_{\phi, \lambda^2\Lambda} = \lambda\cdot \cL_{\phi, \Lambda} \,.
  \]
  Since as $\lambda\searrow 0$,
  \[
    \lambda\gamma_j(\lambda^{-1}s)\to |s| \ \ \text{ if }s<0; \qquad
    \lambda\gamma_j(\lambda^{-1}s)\to |s|e^{i\phi_j} \ \ \text{ if }s>0 \,.
  \]
  Hence the blowdown limit of $[\cL_{\phi, \Lambda}]$ as an integral cycle is $[P_0]+[Q_0]$. Here
  \[
    P_0:= \R^n \oplus 0; \qquad
    Q_0:=\{\Phi t: t\in\R^n\},  \ \text{ where }\Phi:= \diag(e^{i\phi_1},\ldots,e^{i\phi_n})\in \C^{n\times n}.
  \]
  The orientation of $P_0$ is standard, and that of $Q_0$ is the negative of the standard orientation. Both $P_0, Q_0$ are calibrated by $\re\Omega$, since $\sum_j\phi_j=\pi$. 
  Moreover, as $z\to \infty$ on $\cL_{\phi, \Lambda}$, we have $\dist(z, P_0\cup Q_0) = O(|z|^{1-n})$, hence following the terminology in Section \ref{sec:calibration}, $\cL_{\phi, \Lambda}$ has center $0$.
\medskip

  We shall use a different description for the Lawlor necks.
  In the following, fix $\phi, \Lambda$ and denote for simplicity $\cL:=\cL_{\phi, \Lambda}$. Use the unitary map $\mathfrak i(z):=iz$, and give $\mathfrak i(\cL)$ its pushforward orientation.
  Note that the blowdown limit of $\mathfrak i_\#[\cL]$ as an integral cycle is $[P]+[Q]$, here 
  \begin{align}
    P:= 0\oplus i\R^n; \qquad 
    Q:= \{x+iAx: x\in \R^n\}, \ \text{ where } A:=-\diag(\cot\phi_1, \dots, \cot\phi_n) \,;  \label{Equ_Lawlor-pair-planes}
  \end{align}
  with the orientation induced by pushing forward $[P_0], [Q_0]$ by $\mfk i$.  This places the Lawlor neck in the coordinates fixed in Section~\ref{sec:calibration}.
  \begin{proposition}[Global graph of the Lawlor neck]\label{prop:model-graph}
    Under the set-up above, $\mathfrak i(\cL)$ is a global gradient graph over $P^\perp\setminus\{0\}=(\R^n\oplus 0)\setminus \{0\}$, i.e. there exists $U\in C^\infty(\R^n\setminus \{0\}, \R)$ such that
    \begin{equation}\label{eq:model-graph}
      \mathfrak i(\cL) = \{x+i\nabla U(x): x\neq 0\} \,.
    \end{equation}
    In particular, $U$ solves the equation
    \begin{equation}\label{eq:model-phase}
      F(\nabla^2 U(x))=\sum_{j=1}^n\arctan\lambda_j(\nabla^2U(x)) = -\frac{(n-2)\pi}2 \,.
    \end{equation}
  \end{proposition}
  \begin{proof}
    Write $\gamma_j=:a_j+ib_j$, $1\leq j\leq n$. The global graphicality follows from the claim below:
    \begin{Claim} \label{Claim:b_j'>0}
      For each $j$, $b_j'>0$ on $\R$; and $b_j(s)\to 0$ as $s\to -\infty$, $b_j(s)\to +\infty$ as $s\to +\infty$.
    \end{Claim}
    In fact, let $\pi_x: \C^n\to \R^n$, $x+iy\mapsto x$ be the projection onto $x$-factor. Then composing with the parametrization $F_{\phi, \Lambda}$ of $\cL$ and $\mfk i$ gives, \[
      \pi_x\circ \mfk i\circ F_{\phi, \Lambda}(s, u) = (-b_1(s)u_1, \dots, -b_n(s)u_n): \Sph^{n-1} \times \R \to \R^n
    \] 
    Claim \ref{Claim:b_j'>0} guarantees that this is a diffeomorphism onto its image $\R^n\setminus \{0\}$.

    \begin{proof}[Proof of Claim \ref{Claim:b_j'>0}.]
      It follows by definition \eqref{eq:lawlor-functions} that $b_j>0$ and $b_j(-\infty) = 0$, $b_j(+\infty) = +\infty$. Also, a direct calculation yields \[
        a_j'(s) + ib_j'(s) = \gamma_j'(s) = \frac{e^{i\psi_j(s)}}{\sqrt{\rho_j+s^2}}\left(s + i\cP(s)^{-1/2}\right) \,.
      \]
      Hence we need to study the argument $\vartheta(s):=\arg\big( s + i\cP(s)^{-1/2} \big)$. Note that $1+\cot^2\vartheta(s)=1+s^2\cP(s) = F(s)$, taking derivatives on both sides and \eqref{eq:lawlor-polynomial} yield \[
        \vartheta'(s) = -\cP(s)^{-1/2}\cdot \sum_{j=1}^n \frac1{\rho_j+s^2}
        = -\left( \sum_{j=1}^n \psi_j(s) \right)' \,.
      \]
      Since $\vartheta(-\infty)=\pi$ and $\psi_j(-\infty)=0$, we have $\vartheta + \sum_j \psi_j \equiv \pi$. And since $\psi_j(s)\in (0, \phi_j)$, there holds \[
        \arg(\gamma_j'(s)) = \psi_j(s) + \vartheta(s) 
        = \pi-\sum_{k\ne j}\psi_k(s) \in (0, \pi)\,.
      \]
      This proves $b_j'(s)>0$ and thus the Claim.
    \end{proof}

    We have seen that $\mathfrak i(\cL)$ is a global graph $\{y=y(x)\}$. Since $\cL$ is Lagrangian and $\mathfrak i\in U(n)$, the one-form $\sum_j y_j\,dx_j$ is closed. Hence $y=\nabla U$ for some $U\in C^\infty(\R^n\setminus \{0\})$ since $\R^n\setminus\{0\}$ is simply connected for $n\geq 3$. 
    
    Recall that $\mathfrak i_\#[\cL]$ and its two asymptotic planes are calibrated by $\re(i^{-n}\Omega)$. In particular, $U$ solves the special Lagrangian equation \[
      \sum_{j=1}^n \arctan \lambda_j(\nabla^2 U(x)) \equiv \text{constant} \,.
    \]
    And since $\mathfrak i_\#[\cL]$ has blow-down limit $\mfk i_\#([P_0]+[Q_0]) = [P]+[Q]$, we have $\lambda\nabla U(\lambda^{-1} x) \to Ax$ in $C^\infty_{\loc}(\R^n\setminus \{0\})$ as $\lambda\searrow 0$. Thus by \eqref{Equ_Lawlor-pair-planes}, \[
      \sum_{j=1}^n \arctan \lambda_j(\nabla^2 U(x)) 
      \equiv \sum_{j=1}^n \arctan \lambda_j(A) 
      \equiv \sum_{j=1}^n -\frac{\pi}2 + \phi_j = -\frac{(n-2)\pi}{2} \,.
    \]
    This proves \eqref{eq:model-phase}.
  \end{proof}

  \begin{remark}
    It's worth noting that although $Q$ and $\mfk i(\cL)$ are graphs over $\R^n$, their orientation induced by pushing forward by $\mfk i$ differs from the canonical graph orientation by a multiple of $(-1)^{n-1}$. To see this, just notice that under the graph orientation, the (oriented) phase is $F(A)=F(\nabla^2U) = -(n-2)\pi/2$ mod $2\pi$, while Lagrangian planes calibrated by $\re(i^{-n}\Omega)$ should have (oriented) phase $n\pi/2$ mod $2\pi$.
  \end{remark}

\section{Proof of Theorem~\ref{thm:main}}\label{sec:comparison}
  Let $P, Q, T$ be given by Theorem \ref{thm:main}.
  Section~\ref{sec:dimension-two} has already proved the theorem in dimension $n=2$. In the following, let $n\geq 3$.
  Suppose that $T\neq [P]+[Q]$. Then Lemma \ref{lem:structure}, Corollary~\ref{prop:components} and Proposition \ref{prop:coefficient} show that after a translation $\Reg T$ is connected, $\spt T\cap (0\oplus i\R^n) = \emptyset$, and $T$ is a gradient graph over $P\cup Q$ outside a large ball, with the asymptotic estimates for potentials $g_E$: There exists $\kappa>0$ such that 
  \begin{align}
    g_P(t) = -\kappa |t|^{2-n} + O_\infty(|t|^{1-n})\,, \qquad 
    g_Q(t) = \kappa |t|^{2-n} + O_\infty(|t|^{1-n})\,. \label{Equ_PfMain_g_E-expansion}
  \end{align}
  
  Also, let $\hat\cL:= \mfk i(\cL) \subset \C^n$ be a Lawlor neck from Section \ref{sec:model} also asymptotic to $[P]+[Q]$ near infinity. After possibly a rescaling, we may assume that outside a larger ball, $\hat\cL$ is also a gradient graph over $P\cup Q$ with potential $\hat g_E$ satisfying the same asymptotic estimates 
  \begin{align}
    \hat g_P(t) = -\kappa |t|^{2-n} + O_\infty(|t|^{1-n})\,, \qquad 
    \hat g_Q(t) = \kappa |t|^{2-n} + O_\infty(|t|^{1-n})\,.  \label{Equ_PfMain_hat g_E-expansion}
  \end{align}
  
  Let $U\in C^\infty(\R^n\setminus \{0\})$ be the graphical potential of $\hat\cL$ over $(\R^n\oplus 0)\setminus \{0\}$ from Proposition \ref{prop:model-graph}. For $(x, y)\in (\R^n\setminus \{0\})\oplus \R^n$, define 
  \begin{equation} \label{eq:comparison-function}
    \hat q(x,y):= \langle x, y \rangle - 2U(x) \,.
  \end{equation}
  Then by the special Lagrangian equation \ref{prop:model-graph} \eqref{eq:model-phase} and the trace inequality Lemma \ref{Lem_Tr-ineq}, 
  \begin{equation}\label{eq:comparison-defect}
    \Delta_T \hat q(z) = \tr_{T} \nabla^2\hat q(z) = -\tr_{T_z T}M_{\nabla^2 U(x)} \geq 0,  \qquad \forall\, z=(x, y)\in \Reg T \,; 
  \end{equation}
  and equality holds at $z=(x, y)$ if and only if $T_z T\in \{P, T_{x+i\nabla U(x)}\hat\cL\}$ as an unoriented plane. 
  \begin{Claim} \label{Claim_int-Delta-H=0}
    There holds \[
      \int_{\R^{2n}} \tr_T \nabla^2\hat q\ d\|T\| = 0 \,.
    \]
  \end{Claim}
  Note that with this Claim, we immediately know that at every $z=(x, y)\in \Reg T$, $T_z T\in \{P, T_{x+i\nabla U(x)}\hat\cL\}$ as an unoriented plane. Since $\Reg T$ is connected and $T$ can't be a single plane $P$, we must have 
  \begin{align*}
    T_zT=\{(t,D^2U(x)t):t\in \R^n\}\qquad \forall\, z=(x,y)\in \Reg T \,.
  \end{align*} 
  Thus $\nabla (y - DU(x))$ vanishes on $\Reg T$, and then $y - DU(x) \equiv \text{constant}$ on $\Reg T$.
  Hence possibly after a translation, $\Reg T\subset \hat\cL$. By the constancy theorem of stationary integral varifolds, the connectivity of $\hat\cL$, and that $T$ is asymptotic to $[P]+[Q]$ at infinity, we eventually conclude that $T$ is a translation of $[\hat\cL]$, which proves Theorem \ref{thm:main}.

  \begin{proof}[Proof of Claim \ref{Claim_int-Delta-H=0}.]
    Let $\hat Z_E$ be the graphical parametrization of $\hat\cL$ as in Lemma \ref{lem:structure} with $g_E$ replaced by $\hat g_E$. We first prove that, on each end $E\in\{P,Q\}$,
    \begin{equation} \label{eq:end-H-comparison}
      |\nabla_t(\hat q\circ\hat Z_E)(t)|=O(|t|^{1-n}),\qquad
      |\nabla_t(\hat q\circ Z_E-\hat q\circ\hat Z_E)(t)|=O(|t|^{-n}).
    \end{equation}
    
    To prove the first inequality, fix an end $E$ and denote the coordinate components of $Z_E$ and $\hat Z_E \in \R^n \oplus i\R^n$ by $(\bx, \by)$ and $(\hat \bx,\hat \by)$ respectively (viewed as functions on the end of $E$). Then 
    \begin{align*}
      d(\hat q\circ\hat Z_E)
      = \langle d\hat \bx, \hat \by\rangle + \langle\hat \bx, d\hat \by\rangle - 2\langle \nabla U(\hat \bx), d\hat \bx\rangle
      & = \langle\hat \bx, d\hat \by\rangle-\langle\hat \by,d\hat \bx\rangle \\
      & = \langle \nabla^2\hat g_E(t)[t, \cdot] - \nabla\hat g_E(t),\ dt\rangle \,.
    \end{align*}
    here the second equality holds since $\hat\cL$ is the graph of $\nabla U$, and thus $\hat \by=\nabla U(\hat \bx)$; 
    the last equality follows from the unitary changes of coordinates. Since by \eqref{Equ_PfMain_hat g_E-expansion} $D\hat g_E=O(|t|^{1-n})$ and $D^2\hat g_E=O(|t|^{-n})$, this proves the first estimate in \eqref{eq:end-H-comparison}.
    
    Similarly,
    \begin{align*}
      d(\hat q\circ Z_E-\hat q\circ\hat Z_E) 
      & = (\langle d\bx, \by\rangle + \langle \bx, d\by\rangle - 2\langle \nabla U(\bx), d\bx\rangle) 
      - (\langle\hat \bx, d\hat \by\rangle - \langle\hat \by, d\hat \bx\rangle ) \\ 
      & = (\langle \bx, d\by\rangle - \langle \by, d\bx\rangle)
          - (\langle\hat \bx, d\hat \by\rangle - \langle\hat \by, d\hat \bx \rangle)
          + 2\langle \by-\nabla U(\bx), d\bx\rangle \,.\\
      & = \langle \nabla^2(g_E - \hat g_E)|_t[t, \cdot] - \nabla(g_E-\hat g_E)(t),\ dt\rangle + 2\langle \by - \nabla U(\bx), d\bx\rangle \,. 
    \end{align*}
    The first term in the last line is $O(|t|^{-n})$ by \eqref{Equ_PfMain_g_E-expansion} and \eqref{Equ_PfMain_hat g_E-expansion}. It therefore suffices to prove
    \begin{align}
      |\by-\nabla U(\bx)|\cdot |\nabla_t \bx|=O(|t|^{-n}) \qquad \text{ as $t\to \infty$ in $E$, }\ E\in \{P, Q\}\,. \label{Equ_PfMain_|y-U(x)|-asymp}
    \end{align}

    On the $P$ end, $\by(t) = t$, and the end expansions \eqref{Equ_PfMain_g_E-expansion} give
    \begin{align}
      \bx(t) = (2-n)\kappa|t|^{-n} t + O(|t|^{-n}), \qquad 
      \nabla_t \bx(t) = O(|t|^{-n}) \,. \label{Equ_PfMain_P-end_x-asymp}
    \end{align}
    Put $v:=\nabla U(\bx(t))$. Then $(\bx(t), v)\in \hat\cL$ lies on the $P$ end for large $|t|$, with $|v|\to\infty$. 
    Using its $y$-coordinate $v$ as the base parameter gives
    \[
      \bx(t) = (2-n)\kappa |v|^{-n}v + O(|v|^{-n}).
    \]
    Comparing norms first gives $|v|\sim |t|$, and then
    \[
      \big||v|^{1-n} - |t|^{1-n}\big| = O(|t|^{-n}),\qquad
      \big||v| - |t|\big| = O(1).
    \]
    Further comparing directions gives $v-t= O(1)$, that is $|\nabla U(\bx) - \by|= O(1)$ as $t\to \infty$ in $P$, which, together with \eqref{Equ_PfMain_P-end_x-asymp}, proves \eqref{Equ_PfMain_|y-U(x)|-asymp} on the $P$ end.

    On the $Q$ end, let $\pi_x, \pi_y: \R^{2n}\to \R^n$ be the orthogonal projection onto the real and imaginary components. Then for $t\in Q$ tending to infinity, \[
      \hat\bx(t) = \pi_x(t) + O_\infty(|t|^{1-n}) \to \infty, \qquad
      \hat\by(t) = \pi_y(t) + O_\infty(|t|^{1-n}), 
    \] 
    and from \eqref{Equ_PfMain_g_E-expansion} and \eqref{Equ_PfMain_hat g_E-expansion}, \[
      |\hat \bx - \bx| + |\hat \by - \by| = O(|t|^{-n}) \,.
    \]
    Also we know that $\nabla U(\hat\bx) = \hat \by$, and when $x\to \infty$, $\nabla^2 U(x) = A + o(1) = O(1)$, thus \[
      \by - \nabla U(\bx) = (\by-\hat \by) + (\nabla U(\hat \bx)-\nabla U(\bx)) = O(|t|^{-n}) \,.
    \]
    Combine this with $\nabla_t\bx = \pi_x|_Q + o(1) = O(1)$ proves \eqref{Equ_PfMain_|y-U(x)|-asymp} on the $Q$ end as well, and finishes the proof of the second inequality in \eqref{eq:end-H-comparison}.
\medskip


    Note that by Proposition \ref{prop:coefficient}, $\spt T$ is disjoint from $\{x =0\}$, so $\hat q$ is smooth on a neighborhood of it. 
    For $R\gg 1$, by the first variation formula applied to $T_R$ defined in Remark \ref{Rem_T_R,conormal}, 
    \begin{align*}
      \int_{T_R} \tr_T \nabla^2\hat q\ d\|T\| 
      & = \sum_{E\in \{P, Q\}}\int_{\partial_E T_R}\partial_{\nu_R} \hat q\ d\cH^{n-1} \\
      & = \sum_{E\in \{P, Q\}} \int_{\partial_E \hat\cL_R}\partial_{\hat\nu_R} \hat q\ d\cH^{n-1}\  + o(1)
      = \int_{\hat\cL_R} \tr_{\hat\cL} \nabla^2\hat q\ d\|\hat\cL\| + o(1) \,;
    \end{align*}
    where $\nu_R$ and $\hat\nu_R$ denote the outward unit normal of $\partial T_R\subset \spt T$ and $\partial \hat\cL_R \subset \hat \cL$ respectively. 
    Notice that $\tr_{\hat\cL}\nabla^2\hat q = 0$ along $\hat\cL$, sending $R\to +\infty$ proves the Claim.
  \end{proof}

\bibliographystyle{plain}
\bibliography{GMT1103}

\end{document}